\documentclass[10pt]{article}
\usepackage{amssymb,amsmath,amsthm,latexsym, bbm,thmtools}
\usepackage[colorlinks=true, linkcolor=blue, citecolor=blue]{hyperref}
\usepackage[nameinlink]{cleveref}
\usepackage{enumerate}
\usepackage{color}
\usepackage{graphicx}
\usepackage{mathrsfs}
\usepackage[mathscr]{euscript}
\usepackage{amscd}
\allowdisplaybreaks

\newtheorem{thm}{Theorem}[section]

\newtheorem{lemma}[thm]{Lemma}
\newtheorem{prop}[thm]{Proposition}
\newtheorem{defn}[thm]{Definition}

\numberwithin{equation}{section}

\def\Ex {{\mathbb E}}

\def\R {{\mathbb R}}

\def\N {{\mathbb N}}

\def\F {{\mathcal F}}
\def\P {{\mathbb P}}

\def\T {{\mathcal T}}

\def\H {{\mathcal H}}

\def\V {{\mathcal V}}

\newcommand{\norm}[1]{\lVert#1\rVert}
\newcommand{\qv}[1]{\langle#1\rangle}
\def\wt{\widetilde}
\def\wh{\widehat}

\def\op{\mathcal L}
\def\1{{\mathbf 1}}

\def\<{{\langle}}
\def\>{{\rangle}}
\def\eps{\varepsilon}
\def\wh{\widehat}

\def\proof{{\medskip\noindent {\bf Proof. }}}

\def\qed{{\hfill $\square$ \bigskip}}

\makeatletter
\def\tank#1{\protected@xdef\@thanks{\@thanks
		\protect\footnotetext[0]{#1}}}
\def\bigfoot{
	
	\@footnotetext}
\makeatother

\begin{document}
		\title
{\textbf{ Exponential stability of 3D stochastic primitive equations driven by fractional noise} }
\author{
	Lidan Wang\thanks{School of Statistics and Data Science, Nankai University, Tianjin, 300071,  P.R. China. }
	\tank{E-mail:lidanw.math@gmail.com.}
	\qquad
	Guoli Zhou\thanks{Corresponding author at School of Statistics and Mathematics, Chongqing University, Chongqing, 400044, P.R. China. }
	\tank{E-mail:zhouguoli736@126.com.}
}
\date{}
\maketitle
\begin{abstract}
	In this article, we study the stability of solutions to 3D stochastic primitive equations driven by fractional noise. Since the fractional Brownian motion is essentially different from Brownian motion, lots of stochastic analysis tools are not available to study the exponential stability for the stochastic systems. Therefore, apart from the standard method for the case of Brownian motion, we develop a new method to show that 3D stochastic primitive equations driven by fractional noise converge almost surely exponentially to the stationary solutions. This method may be applied to other stochastic hydrodynamic equations and other noises including Brownian motion and L\'evy noise. 
\end{abstract}	

\bigskip


\bigskip
	
%%%%%%%%%%%%%%%%%%%%%%
\section{Introduction}
%%%%%%%%%%%%%%%%%%%%%%	
In this article, we study the exponential stability of solutions to 3D stochastic primitive equations (SPEs) driven by fractional Brownian motion. We first define  a cylindrical domain $\mho=M\times(-h,0)\subset\R^3$, where $M\subset\R^2$ is a smooth bounded domain, then formulate 3D stochastic PEs of Geophysical Fluid Dynamics as follows:
\begin{align*}
	&\partial_t v+L_1v+(v\cdot \nabla)v+w\partial_zv+fv^{\perp}+\nabla p=G_1(t)\dot W_1^H,\\
	&\partial_tT+L_2T+(v\cdot \nabla) T+w\partial_z T=Q_2+G_2(t)\dot W_2^H,\\
	&\nabla\cdot v+\partial_zw=0,\\
	&\partial_zp+T=0.
\end{align*}
The unknowns are the fluid velocity field $(v,w)=(v_1,v_2,w)\in\R^3$, with $v=(v_1,v_2), v^\perp=(-v_2,v_1)$ being horizontal, the temperature $T$ and the pressure function $p$.

The Coriolis parameter $f$ is defined by $f=f_0+\beta y$ and $Q_2$ is a given heat sourse, the viscocity and the heat diffusion operators $L_1,L_2$ are given by
$$L_i=-\nu_i\Delta-\mu_i\partial_{zz},\ i=1,2.$$
Here the positive constants $\nu_1,\mu_1$ are the horizontal and vertical Reynolds numbers, and positive constants $\nu_2,\mu_2$ stand for the horizontal and vertical heat diffusivity.

Throughout this paper, $\nabla,\Delta,\text{div}$ represent the horizontal gradient, Laplacian and divergence, respectively. $\dot W_i^H(t,x,y,z), i=1,2$ stands for the informal derivative for the fractional Wiener process $W_i^H$ that will be introduced later in \autoref{subsection:fBM}.

For the cylindrical domain, the boundary can be partitioned into three parts: $\partial\mho:=\Gamma_u\cup\Gamma_b\cup\Gamma_s$, where
\begin{align*}
	&\Gamma_u=\{(x,y,z)\in\bar\mho:z=0\},\\
	&\Gamma_b=\{(x,y,z)\in\bar\mho:z=-h\},\\
	&\Gamma_s=\{(x,y,z)\in\bar\mho:(x,y)\in\partial M,-h\leq z\leq 0\}.
\end{align*}  	
	We now supplement the SPE model with initial and boundary conditions:
	\begin{align*}
		&v(x,y,z,0)=v_0(x,y,z),\ T(x,y,z,0)=T_0(x,y,z),\\
		&\partial_zv=\eta,w=0,\partial_zT=-\alpha(T-\tau)\text{ on }\Gamma_u,\\
		&\partial_zv=0,w=0,\partial_zT=0\text{ on }\Gamma_b,\\
		&v\cdot\vec{n}=0, \partial_{\vec{n}}v\times\vec{n}=0,\partial_{\vec{n}}T=0\text{ on }\Gamma_s,
	\end{align*}
	where $\eta(x,y)$ represents the wind stres on the surface of the ocean, $\alpha>0$, $\tau$ is the typical temperature distribution on the ocean surface, and $\vec{n}$ is the norm vector on $\Gamma_s$. For the sake of simplicity, we assume $Q_2$ is independent of time, and $\eta=\tau=0$. We would like to mention that the results presented in this paper can be extended to the general case provided with simple modifications.
	
	By easy calculations, we have that
	\begin{align*}
		&w(x,y,z,t)=-\int_{-1}^z\nabla\cdot v(x,y,\xi,t)d\xi,\\
		&p(x,y,z,t)=p_s(x,y,z,t)-\int_{-1}^zT(x,y,\xi,t)d\xi.
	\end{align*}
	Hence, the stochastic model, together with initial and boundary conditions, can be reformulated into the following equivalent form:
	\begin{align}
	\label{eqn:1.1}
		&\partial_t v+L_1v+(v\cdot \nabla)v-\left(\int_{-1}^z\nabla\cdot v(x,y,\xi,t)d\xi\right)\partial_zv\notag\\
		&\hspace{1cm}+fv^{\perp}+\nabla p_s(x,y,z,t)-\int_{-1}^z\nabla T(x,y,\xi,t)d\xi=G_1(t)\dot W_1^H,\\
		&\partial_tT+L_2T+(v\cdot \nabla )T-\left(\int_{-1}^z\nabla\cdot v(x,y,\xi,t)d\xi\right)\partial_z T=Q_2+G_2(t)\dot W_2^H,\\
		\label{eqn:1.3}
		&\partial_zv|_{\Gamma_u}=\partial_zv|_{\Gamma_b}=0,\ v\cdot\vec{n}|_{\Gamma_s}=0,\ \partial_{\vec{n}}v\times\vec{n}|_{\Gamma_s}=0,\\
		\label{eqn:1.4}
		&(\partial_zT+\alpha T)|_{\Gamma_u}=\partial_zT|_{\Gamma_b}=0,\ \partial_{\vec{n}}T|_{\Gamma_s}=0,\\
		\label{eqn:1.5}
		&v(x,y,z,0)=v_0(x,y,z),\ T(x,y,z,0)=T_0(x,y,z).
	\end{align}
The Primitive equations are the basic model used in the study of climate and
weather prediction, which describe the motion of the atmosphere when the hydrostatic
assumption is enforced $\cite{G,H1,H2}$. As far as we know, their mathematical study was initiated by J. L. Lions, R. Teman and S. Wang$(\cite{LTW1}-\cite{LTW4})$. And this research field has developed and  has received
considerable attention from the mathematical community over the last two decades. Taking advantage of the fact that the pressure is essentially two-dimensional in the PEs, Cao and Titi $\cite{CT1}$ proved the global results for the existence of strong solutions of the full
three-dimensional PEs. Independently, I. Kukavica and M. Ziane $\cite{KZ}$ developed a different proof which allows one to treat non-rectangular domains as well as different, physically realistic, boundary conditions. The existence of
the global attractor is given by Ju $\cite{Ju}$.  For the PEs with partial dissipation, we refer the reader to the
papers $\cite{CIN, CLT1, CLT2, CLT3, CT2}$. There are also other good works about the global well-posedness theory of PEs, we do not intend to mention each of them here.

The breakthroughs in the deterministic case motivated the development of the theory for the stochastic PEs. B. Ewald, M. Petcu, R. Teman $\cite{EPT}$ and N. Glatt-Holtz, M. Ziane $\cite{GHZ}$ considered a two-dimensional stochastic PEs. Then N. Glatt-Holtz and R. Temam $\cite{GHT1, GHT2}$ extended the case to the greater generality of physically relevant boundary conditions and nonlinear multiplicative noise. Following the methods similar to $\cite{CT1}$, Boling Guo and Daiwen Huang $\cite{GH}$ studied the global well-posedness and long-time behavior of the three-dimensional system with  additive noise. Using an approach similar   to $\cite{KZ}$, A. Debussche, N. Glatt-Holtz, R. Temam and M. Ziane \cite{DGTZ} considered three-dimensional system with multiplicative noise. In a subsequent paper, N. Glatt-Holtz, I. Kukavica, V. Vicol, and M. Ziane \cite{GKVZ} established the existence and regularity of invariant measure for the SPEs. The ergodic theory and large deviations for the 3D SPEs were obtained by Dong, Zhai and Zhang
in \cite{DZZ1,DZZ2}. Some analytical properties of weak solutions of 3D stochastic primitive equations with
periodic boundary conditions were obtained in \cite{DZr}, in which the martingale problem associated to
this model is shown to have a family of solutions satisfying the Markov property. Concerning the exponential behavior of SPEs, T. T. Medjo established the stability result for SPEs in \cite{M1} when the noise is Brownian motion. H. J. Gao and C. F. Sun also studied the long-time behavior,  asymptotic and regular properties of SPEs, see \cite{GS1, GS2, GS3}.
\par
As it is pointed out in \cite{P, MTV}that studies from climate show that the complex multi-scale nature of the earth's climate system results in many uncertainties that should be accounted for in the basic dynamical models of atmospheric and oceanic processes. It is further suggested in \cite{Zgl} that the  uncertainties prefer to be $\mathbf{Non}$-$\mathbf{Markovian}$. This is the motivation for us to study the stability of SPEs driven by fractional Brownian motion.
 \par
 Since the noise here is fractional Brownian motion which has memory for the increment and is essentially different from Brownian motion, the stochastic integral is not an It\^{o} integral. Therefore, we instead define stochastic integrals via pathwise generalized Stieltjes integrals as is the case in $\cite{MN,NR,Z}$, for the details, one can check  Section \ref{sec:2} in this article. As a result, the method studying the stability of the stationary solutions to SPEs in the present work is different from the standard method, see \cite{ CLTa1, CLTa2, M1} and other references. In the following, we will illustrate the differences more clearly.

 \par
 As we know if one tries to obtain the moment stability for the stochastic equations with nonlinear multiplicative Wiener noise,  It\^{o} formula will play an important role, see  Theorem 3.2 in \cite{M1} and other references. Furthermore, if one  tries to establish stability of the sample paths, then the moment stability and Borel-Cantelli lemma are key tools. One can refer to Theorem 3.3 in \cite{M1} and other references. If the noise is a spacial linear noise as considered in Theorem 4.2 of \cite{M1}, Medjo takes advantage of the polynomial growth of Brownian motion to show that the stationary solution to SPEs is the almost surely exponentially stable.

 \par
 In our article, we consider the case of additive fractional Brownian motion, as is the case in \cite{CLTa1, CLTa2} and \cite{M}. A common method to deal with this type of noise is to introduce a  fractional Ornstein-Uhlenback process to convert the SPEs into primitive equations with random coefficients. Here, we should point out that under the definition of the stochastic integral the fractional Ornstein-Uhlenback process is uniformly bounded when time $t$ goes to infinite, please see \autoref{prop:ZH3} for details. This is a key result which helps us to show that the stationary solution to SPEs is exponentially stable, which is stated in \autoref{thm:expstability}.

 \par
  It is known that one other specialty about Ornstein-Uhlenback processes driven by Brownian motion is  the polynomial growth of sample paths. Therefore, the uniform Gronwall inequality can not be applied to stochastic equations. For this reason, the random attractor obtained for the stochastic equations are pullback random attractor. There are  many literatures about the topic of random attractor, we list some of them here for the convenience of readers,  \cite{BCL, BLW, CDF, CF, GGS}.  However, in this article, we find that the uniform boundedness of the fractional Ornstein-Uhlenback process may open a way to obtain global attractor (not the pullback random attractor) for the stochastic equations. By virtue of the uniform boundedness of the fractional Ornstein-Uhlenback process, we establish the uniform estimate for the strong solution to SPEs via uniform Gronwall inequality, see \autoref{prop:UH1}. More discussions about the existence of global attractors are given in authors' another recent work \cite{WZ}. 
        
 \par
 To study the long-time behavior of stochastic dynamic system, the ergodicity of fractional Ornstein-Uhlenback process is an important tool which can help us to obtain the uniform \textit{a priori} estimates with respect to time $t$ for the solutions. But the dissipation of SPEs is not enough (see the proof of \autoref{prop:UH1}), one can not  directly use the ergodicity of fractional Ornstein-Uhlenback process to establish the desired energey estimates. Therefore, to overcome the difficulty, we introduce another  fractional Ornstein-Uhlenback process depending on the parameter $\beta$ to obtain uniform estimates,  please see \eqref{eqn:2.9} for the definition of this spacial fractional Ornstein-Uhlenback process and detailed calculations below \eqref{eqn:3.5}.

 \par 
 Based on the foregoing results, we establish our main result \autoref{thm:expstability} which shows that the stationary solution to SPEs is exponentially stable. One of the difficulties is to deal with the stochastic term. As we know  fractional Brownian motion is not a semi-martingale, thus, we can not use the classic tools, i.e.,   It\^{o} formula and Burkholder-Davis-Gundy inequality to achieve the stability results. Our idea is to discretize the stochastic integral with respect to fractional Brownnian motion, which helps us to open the way to make full use of the stationary properties and polynomial growth properties of the increments of the fractional Brownian motion as well as regularities of the fractional Brownian motion. After delicate and careful estimates via the regularities and stationary properties of fractional Brownian motion, we establish the uniform estimate of the stochastic terms, please see the estimates of $K_{1},K_{2}, N_{1}$ and $N_{2}$  in the proof of \autoref{thm:expstability}. 
 
 \par
 We would like to mention that the result and method presented in this article may be a basic tool to study the stability behavior of stochastic partial differential equations driven by fractional Brownian motion.

	The structure of the paper is as follows: In  \autoref{sec:2}, we give definitions of functional spaces and operators regards to the SPEs, and introduce the pathwise integral with fractional calculus techniques, with those techniques, under certain conditions on the forcing terms, we show that the O-U processes are uniformly bounded in $H^3$ norm. In \autoref{sec:3}, we give the definition of exponential stability, then discuss the uniform bounded property for $H^1$ norm of solutions, and finally under a stronger condition for the forcing terms, we establish the exponential stability in the almost sure sense for the solutions of SPEs.
%%%%%%%%%%%%%%%%%%%%%%
\section{Preliminaries }
\label{sec:2}
%%%%%%%%%%%%%%%%%%%%%%

%%%%%%%%%%
\subsection{Functional spaces}
%%%%%%%%%%
For $1\leq p\leq \infty$, let $L^p(\mho), L^p(M)$ be the usual Lebesgue spaces with norms $|\cdot|_{p}$ and $|\cdot|_{L^p(M)}$, respectively. For $m>0$, we denote by $(H^{m,p}(\mho),\norm{\cdot}_{m,p})$ and $(H^{m,p}(M),\norm{\cdot}_{H^{m,p}(M)})$ be the usual Sobolev spaces. When $p=2$, we write them as $(H^{m}(\mho),\norm{\cdot}_{m})$ and $(H^{m}(M),\norm{\cdot}_{H^{m}(M)})$ for short. For simplicity of notations, we sometimes directly write $|\cdot|_p, \norm{\cdot}_m$ as norms in $L^p(M), H^m(M)$, if there is no confusion.

We now define the following function spaces:
\begin{align*}
	&V_1=\{v\in(C^\infty(\mho))^2:\partial_zv|_{\Gamma_u}=\partial_zv|_{\Gamma_b}=0, v\cdot\vec{n}|_{\Gamma_s}=0, \partial_{\vec{n}}v\times\vec{n}|_{\Gamma_s}=0,\int_{-1}^0\nabla\cdot vdz=0\};\\
	&V_2=\{T\in C^\infty(\mho):(\partial_zT+\alpha T)|_{\Gamma_u}=\partial_zT|_{\Gamma_b}=0, \partial_{\vec{n}}T|_{\Gamma_s}=0\}.
\end{align*}
Denote by $\mathcal V_1, \mathcal V_2$ the closure spaces of $V_1, V_2$ under $(H^1(\mho))^2,H^1(\mho)$, respectively. Let $\mathcal H_1$ be the closure space of $V_1$ under the norm $|\cdot|_2$ and $\mathcal H_2=L^2(\mho)$. Now set
$$\mathcal V=\mathcal V_1\times\mathcal V_2,\ \mathcal H=\mathcal H_1\times\mathcal H_2.$$
For $U:=(v,T),\wt U:=(\wt v,\wt T)\in\mathcal V$, we equip $\mathcal V$ with the inner product
\begin{align*}
	&\qv{U,\wt U}_{\mathcal V}:=\qv{v,\wt v}_{\mathcal V_1}+\qv{T,\wt T}_{\mathcal V_2},\\
	&\qv{v,\wt v}_{\mathcal V_1}:=\int_\mho(\nabla v\cdot\nabla\wt v+\partial_z v\cdot\partial_z\wt v)dxdydz,\\
	&\qv{T,\wt T}_{\mathcal V_2}:=\int_\mho(\nabla T\cdot\nabla\wt T+\partial_z T\cdot\partial_z\wt T)dxdydz+\alpha\int_{\Gamma_u}T\wt Td\Gamma_u.
\end{align*}
Consequently, the norm in $\mathcal V$ is defined by $\norm{U}_1:=\qv{U,U}_{\mathcal V}^{1/2}$. Similarly, we define the inner product in $\mathcal H$ by
\begin{align*}
	&\qv{U,\wt U}_{\mathcal H}:=\qv{v,\wt v}+\qv{T,\wt T},\\
	&\qv{v,\wt v}:=\int_\mho v\cdot\wt vdxdydz,\ \qv{T,\wt T}=\int_\mho T\wt Tdxdydz.
\end{align*}
%%%%%%%%%%
\subsection{Operators}
%%%%%%%%%%
 We define the bilinear forms $a:\mathcal V\times\mathcal V\to\R, a_i:\V_i\times\V_i\to\R, i=1,2$, by
\begin{align*}
	&a_1(v,\wt v)=\qv{v,\wt v}_{\V_1}, a_2(T,\wt T)=\qv{T,\wt T}_{\V_2}, a(U,\wt U)=\qv{U,\wt U}_\V.
\end{align*}
Denote by $\V',\V_i'$ the dual spaces of $\V,\V_i$ for $i=1,2$, then the corresponding linear operators $A:\mathcal V\to\V', A_i:\V_i\to\V_i', i=1,2$ can be defined as follows:
\begin{align*}
	&\qv{A_1 v,\wt v}=a_1(v,\wt v), \qv{A_2 T,\wt T}=a_2(T,\wt T), \qv{A U,\wt U}=a(U,\wt U),
\end{align*}
where $U=(v,T),\wt U=(\wt v,\wt T)\in \V$.

Now for $i=1,2$, define $D(A_i):=\{\rho\in\V_i,A_i\rho\in H_i\}$. Since $A_i^{-1}$ is a self-adjoint compact operators in $\H_i$, by the classic spectral theory, we can define the power $A_i^s$ for any $s\in\R$. Then $D(A_i)'=D(A_i^{-1})$ is the dual space of $D(A_i)$ and $\V_i=D(A_i^{1/2}), \V_i'=D(A_i^{1/2})$. Moreover, we have the compact embedding relationship
$$D(A_i)\subset\V_i\subset\H_i\subset\V_i'\subset D(A_i)', \ i=1,2.$$
And
$$\norm{\cdot}_1^2=a_i(\cdot,\cdot)=\qv{A_i\cdot,\cdot}=\qv{A_i^{1/2}\cdot,A_i^{1/2}\cdot},\ i=1,2.$$
Hereafter we denote by $\lambda_1>0$ a constant such that
\begin{equation}
\label{eqn:L2H1}
|U|_2\leq\lambda_1\norm{U}_1,\ U\in\V;\hspace{5mm}\norm{U}_1\leq\lambda_1|AU|_2,\ U\in D(A).
\end{equation}
Now we define the nonlinear operator for $U=(v,T),\wt U=(\wt v,\wt T),\wh U=(\wh v,\wh T)\in \V$,
$$b(U,\wt U,\wh U)=\qv{B(U,\wt U),\wh U}=b_1(v,\wt v,\wh v)+b_2(v,\wt T,\wh T),$$
where
\begin{align*}
	&b_1(v,\wt v,\wh v)=\qv{B_1(v,\wt v),\wh v}:=\int_\mho\left[(v\cdot\nabla)\wt v-\left(\int_{-1}^z\nabla\cdot v(x,y,\xi,t)d\xi\right)\partial_z\wt v\right]\cdot\wh vdxdydz,\\
	&b_2(v,\wt T,\wh T)=\qv{B_2(v,\wt T),\wh T}:=\int_\mho\left[(v\cdot\nabla)\wt T-\left(\int_{-1}^z\nabla\cdot v(x,y,\xi,t)d\xi\right)\partial_z\wt T\right]\cdot\wh Tdxdydz.
\end{align*}
There are some basic properties for the nonlinear operator $B$. The proof can be found in \cite{Ju}.
\begin{lemma}
	\label{lemma:2.1}
	There exists a constant $c_0>0$ that is independent of $U,\wt U,\wh U$ such that
	$$b(U,\wt U,\wt U)=\qv{B(U,\wt U), \wt U}=0,\ \text{ for }U\in \mathcal V, \wt U\in D(A).$$
	$$|b(U,\wt U,\wh U)|\leq c_0\norm{U}_1^{1/2}|AU|_2^{1/2}\norm{\wt U}_1^{1/2}|A\wt U|_2^{1/2}\|\wh U\|_1,\ \text{ for }U,\wt U\in D(A), \wh U\in\V.$$
	$$|b(U,\wt U,\wh U)|\leq c_0\norm{U}_1^{1/2}|AU|_2^{1/2}\norm{\wt U}_1\norm{\wh U}_1^{1/2}|\wh U|_2^{1/2}, \ \text{ for }U\in D(A), \wt U,\wh U\in\V.$$
	$$|b(U,\wt U,\wh U)|\leq c_0\norm{U}_1\norm{\wt U}_1^{1/2}|A\wt U|_2^{1/2}\norm{\wh U}_1^{1/2}|\wh U|_2^{1/2}, \ \text{ for }\wt U\in D(A), U,\wh U\in\V.$$
\end{lemma}
We also define another bilinear operator $R:\V\to\V'$ by
$$r(U,\wt U)=\qv{R(U),\wt U}=\qv{fv^{\perp},\wt v}+\qv{\int_{-1}^z\nabla T,\wt v}, \ \text{ for }U=(v,T),\wt U=(\wt v,\wt T)\in\V.$$
Recall that for $U\in\mathcal V, \wt U\in \H$, there exists $\alpha_0>0$ such that
\begin{equation}
\label{eqn:R}
	|r(U,\wt U)|=|\qv{R(U),\wt U}|\leq \alpha_0\norm{U}_1|\wt U|_2
\end{equation}
Finally, we write $Q:=(0,Q_2)$ as the forcing term. We write $G(t)=(G_1(t),G_2(t)$ as the coefficients for the random term. Hereafter we set
$$\bar\nu_i=\min(\nu_i,\mu_i),\ i=1,2;\ \bar\nu=\min(\bar\nu_1,\bar\nu_2).$$
Thus, we consider the stochastic 3D PEs of the ocean, driven by fractional Brownian motion, written in the following abstract mathematical setting
\begin{equation}
\label{eqn:2.3}
	dU(t)=[-\nu AU(t)-B(U(t))-R(U(t))+Q(t)]dt+G(t)dW^H(t),\ U(0)=U_0:=(v_0,T_0).
\end{equation}
%%%%%%%%%%%%
\subsection{Fractional Brownian motion}
\label{subsection:fBM}
%%%%%%%%%%%%
Set $\alpha\in(0,1)$, for a function $f:[0,T]\to\R$ that is regular enough, we define the Weyl fractional derivatives as follows:
\begin{align*}
	D_{0+}^\alpha f(t)=&\frac{1}{\Gamma(1-\alpha)}\left(\frac{f(t)}{t^\alpha}+\alpha\int_0^t\frac{f(t)-f(u)}{(t-u)^{\alpha+1}}du\right),\\
	D_{T-}^\alpha f(t)=&\frac{(-1)^{\alpha}}{\Gamma(1-\alpha)}\left(\frac{f(t)}{(T-t)^\alpha}+\alpha\int_t^T\frac{f(t)-f(u)}{(u-t)^{\alpha+1}}du\right),
\end{align*}
provided the singular integrals on the right hand side exist for almost all $t\in[0,T]$, and $\Gamma$ stands for the Gamma function.

For $\phi\in L^1([0,T];\R)$, we define the left and right hand side fractional Riemann-Liouville integrals of $\phi$ of order $\alpha$ for almost all $t\in(0,T)$ by
\begin{align*}
	I_{0+}^\alpha\phi(t)=&\frac{1}{\Gamma(\alpha)}\int_0^t(t-u)^{\alpha-1}\phi(u)du,\\
	I_{T-}^\alpha\phi(t)=&\frac{(-1)^{-\alpha}}{\Gamma(\alpha)}\int_t^T(u-t)^{\alpha-1}\phi(u)du.
\end{align*}
If $f=I_{0+}^\alpha\phi$, then the Weyl left-side derivative of $f$ exists and $D_{0+}^\alpha f=\phi$. A similar result holds for the right-side fractional integral. See \cite{SKM} for a comprehensive introduction for  the theory of fractional integrals and derivatives.

We now define  $W^{\alpha,1}([0,T];\R)$ as the space of measurable functions $f:[0,T]\to\R$ with
$$\norm{f}_{\alpha,1}:=\int_0^T\left(\frac{|f(s)|}{s^\alpha}+\int_0^s\frac{|f(s)-f(u)|}{(s-u)^{\alpha+1}}du\right)ds<\infty,$$
where $\alpha\in(0,\frac{1}{2})$.

Following \cite{Z}, we define the generalized Stieltjes integral $\int_0^Tfdg$ by
\begin{equation}
\label{eqn:2.4}
	\int_0^Tfdg=(-1)^\alpha\int_0^TD_{0+}^\alpha f(s)D_{T-}^{1-\alpha}g_{T-}(s)ds,
\end{equation}
where $g_{T-}(s)=g(s)-g(T)$. Under the above hypotheses, the above integral exists for all $t\in[0,T]$, and by \cite{NR}, we have
$$\int_0^tfdg=\int_0^Tf1_{(0,t)}dg.$$
Furthermore, we have the following estimate
\begin{equation}
\label{eqn:2.5}
	\bigg|\int_0^tfdg\bigg|\leq C_\alpha(g)\norm{f}_{\alpha,1},
\end{equation}
where
$$C_\alpha(g):=\frac{1}{\Gamma(\alpha)\Gamma(1-\alpha)}\sup_{0<s<t<T}\left(\frac{|g(t)-g(s)|}{(t-s)^{1-\alpha}}+\int_s^t\frac{|g(u)-g(s)|}{(u-s)^{2-\alpha}}du\right).$$
On a complete probability space $(\Omega,\F,\P)$, for $H\in(0,1)$, we let $B^H_i=(B^H_i(t)_{t\in\R})$ be a sequence of independent, identically distributed continuous centered Gaussian process with covariance function
$$R(s,t)=\frac{1}{2}(|t|^{2H}+|s|^{2H}-|t-s|^{2H}),\ s,t\in\R. $$
This type of process is called a two-sided one-dimensional fractional Brownian motion (fBm) with Hurst parameter $H$. When $H=\frac{1}{2}$, $B^{1/2}_i$ is the standard Brownian motion.

For $j\in\{1,2\}$, let $H_j$ be linear, self-adjoint, positive trace-class operators on $\H_j$, that is, for a complete orthonormal basis $(e_{i,j})_{i\in\N_+}$ in $\H_j$, there exists a sequence of nonnegative values $(\lambda_{i,j})_{i\in\N_+}$ such that $\text{tr}H_j=\sum_{i=1}^{\infty}\lambda_{i,j}<\infty$. We now introduce $\H_j$-valued fractional Brownian motion $W_j^H, j=1,2$, with covariance operator $H_j$, and Hurst parameter $H$, as follows:
\begin{equation}
\label{eqn:2.6}
	W_j^H(t):=\sum_{i=1}^{\infty}\sqrt{\lambda_{i,j}}e_{i,j}B_i^H(t).
\end{equation}
For $H>\frac{1}{2}$, take a parameter $\alpha\in(1-H,\frac{1}{2})$ which will be fixed throughout the paper. Let $f\in W^{\alpha,1}([0,T];\R)$, we define $\int_0^Tf(s)dB_i^H(s)$ in the sense of \eqref{eqn:2.4} pathwise. By \cite{NR}, $C_\alpha(B_i^H)<\infty$, $\P$-a.s. for $t\in\R_+, i\in\N_+$.

For $j=1,2$, denote by $\op(\V_j)$ the space of linear bounded operators on $\V_j$, now suppose $l:\Omega\times[0,T]\to\op(\V_j)$  is an operator-valued function such that $le_{i,j}\in W^{\alpha,1}([0,T];\V_j)$ for any $i\in\N_+,\omega\in\Omega$. Now we define
\begin{align}
	\int_0^T\varphi(s)dW_j^H(s):=&\sum_{i=1}^\infty\int_0^Tl(s)H_j^{1/2}e_{i,j}dB_i^H(s)\notag\\
	=&\sum_{i=1}^\infty\sqrt{\lambda_{i,j}}\int_0^Tl(s)e_{i,j}dB_i^H(s),
\end{align}
where the convergence of the above series is understood as $\P$-a.s. convergence in $\V_j$.

%%%%%%
\subsection{Regularity of O-U processes}
%%%%%

Now we consider the stochastic equations for $j=1,2$, $\beta>1$,
\begin{equation}
	dZ_j(t)=(-A_jZ_j-\beta Z_j)dt+G_j(t)dW_j^H(t),\ Z_j(0)=0.
\end{equation}
The solution can be interpreted pathwisely in the mild sense, that is, the solution $(Z_j(t))_{t\in[0,\T]}$ is a $\V_j$-valued process whose paths are elements of the space $W^{\alpha,1}([0,\T];\V_j)$ with probability one, for $\alpha\in(1-H,\frac{1}{2})$, such that,
\begin{equation}
\label{eqn:2.9}
	Z_j(t)=\int_0^te^{-(t-s)(A_j+\beta)}G_j(s)dW_j^H(s).
\end{equation}
We first have the following result which gives the growth rates of $Z_j(t), j=1,2$, under certain conditions on the forcing terms.
\begin{prop}
	\label{prop:ZH3}
	For $j=1,2$, denote by $0<\gamma_{1,j}\leq\gamma_{2,j}\leq\cdots$ the eigenvalues of $A_j$ with corresponding eigenvectors $e_{1,j},e_{2,j},\cdots$, assume the following condition holds
	\begin{equation}
	\label{eqn:2.10}
		\sum_{i=1}^\infty \lambda_{i,j}^{1/2}\gamma_{i,j}^{5/2}<\infty,\ \text{ for }j=1,2.
	\end{equation}
Then for $\T>0$, $(Z_1(t))_{t\in[0,\T]},(Z_2(t))_{t\in[0,\T]}$ exist as generalized Stieltjes integrals in the sense of \cite{Z}, and 
$$(Z_{1}(t))_{t\in [0,\T]}\in C([0,\T]; (H^{3}(\mho))^{2}), a.s.\ \mathrm{and}\   (Z_{2}(t))_{t\in [0,\T]}\in C([0,\T]; H^{3}(\mho)),\ a.s. . $$
	Assume furtherly that the forcing terms $G_1(t), G_2(t)$ only depend on $t$ and satisfy
	\begin{equation}
	\label{eqn:Gcon}
	|G_1(t)|+|G_1'(t)|\leq M_1(1+t)^{-2},\  |G_2(t)|+|G_2'(t)|\leq M_2(1+t)^{-2}.
	\end{equation}
Then for $\beta>0$, $Z_j(t), j=1,2$ are uniformly bounded in $H^3$ norms, in the sense that there exist random variables $C_1(\omega), C_2(\omega)$ taking finite values such that
	$$\sup\limits_{t\in [0,\infty)}\|Z_1(t)\|_{3}\leq C_1(\omega)<\infty, a.s.,\ \ \ \  \sup\limits_{t\in [0,\infty)}\|Z_2(t)\|_{3}\leq C_2(\omega), a.s..$$
\end{prop}
\proof For the proof of the continuity of $Z_{1}$ and $Z_{2},$ one can refer to Propositions 2.1 of \cite{Zgl}. Here we will give a short proof of the uniform boundedness of $Z_{1}$ and $Z_{2}.$ By the definition in \eqref{eqn:2.9} and the property of Riemann-Stieltjes integral (see Theorem 4.2.1 in \cite{Z}), we have
\begin{align*}
	\|Z_1(t)\|_{3}=&\bigg\|\int_0^te^{-(t-s)(A_1+\beta)}G_1(s)dW_1^H(s)\bigg\|_{3}\\
	=&\bigg\|\sum_{k=1}^{\infty}\sqrt{\lambda_{k,1}}\int_0^te^{-(t-s)(\gamma_{k,1}+\beta)}e_{k,1}G_1(s)dB_k^H(s)\bigg\|_{3}\\
	=&\bigg\|\sum_{k=1}^{\infty}\sqrt{\lambda_{k,1}}e_{k,1}\Big[G_1(t)B_k^H(t)\\
	&\hspace{2.5cm}-\int_0^tB_k^H(s)[(\gamma_{k,1}+\beta)e^{-(t-s)(\gamma_{k,1}+\beta)}G_1(s)+e^{-(t-s)(\gamma_{k,1}+\beta)}G'_1(s)]ds\Big]\bigg\|_{3}.
\end{align*}
Since $|B_k^H(t)|\leq t^2+c(\omega)$ (see Lemma 2.6 of \cite{MS}), $|G_1(t)|+|G_1'(t)|\leq M_1(1+t)^{-2}$,  we get
\begin{align*}
	\|Z_1(t)\|_{3}\leq&M_1(1+t)^{-2}(t^2+c(\omega))\sum_{k=1}^{\infty}\sqrt{\lambda_{k,1}}\gamma_{k,1}^{\frac{3}{2}}\\
	&+M_1\sum_{k=1}^{\infty}\sqrt{\lambda_{k,1}}\gamma_{k,1}\gamma_{k,1}^{\frac{3}{2}}\int_0^t(s^2+c(\omega))(1+s)^{-2}	e^{-(t-s)(\gamma_{k,1}+\beta)}ds\\
	\leq &M_1c\sum_{k=1}^{\infty}\sqrt{\lambda_{k,1}}\gamma_{k,1}^{\frac{3}{2}}+M_1c\sum_{k=1}^{\infty}\sqrt{\lambda_{k,1}}\gamma_{k,1}^{\frac{5}{2}}\int_0^t	e^{-s(\gamma_{k,1}+\beta)}ds\\
	\leq&C_1(\omega).
	\end{align*}
The estimate for $Z_2(t)$ follows similarly.

\qed
%%%%%%%%%%%%%%
\subsection{Global well-posedness}
%%%%%%%%%%%%%
Let $u(t)=v(t)-Z_1(t)$ and $\theta(t)=T(t)-Z_2(t)$, a stochastic process $U(t,\omega)=(v,T)$ is a strong solution to \eqref{eqn:1.1}-\eqref{eqn:1.5} on $[0,\T]$, if and only if $(u,\theta)$ is a strong solution to the following problem on $[0,\T]$:
\begin{align}
\label{eqn:2.12}
	&\partial_tu+L_1u+[(u+Z_1)\cdot\nabla](u+Z_1)+w(u+Z_1)\partial_z(u+Z_1)\notag\\
	&+f(u+Z_1)^\perp+\nabla p_s-\int_{-1}^z\nabla (\theta+Z_2)(x,y,\xi,t)d\xi=\beta Z_1;\\
	&\partial_t\theta+L_2\theta+[(u+Z_1)\cdot\nabla](\theta+Z_2)+w(u+Z_1)\partial_z(\theta+Z_2)=Q+\beta Z_2;\\
	&\int_{-1}^0\nabla\cdot udz=0;\\
	&\partial_zu|_{\Gamma_u}=\partial_zu|_{\Gamma_b}=0,\ u\cdot\vec{n}|_{\Gamma_s}=0,\ \partial_{\vec{n}}u\times\vec{n}|_{\Gamma_s}=0;\\
	&(\partial_z\theta+\alpha\theta)|_{\Gamma_u}=\partial_z\theta|_{\Gamma_b}=0,\ \partial_{\vec{n}}\theta|_{\Gamma_s}=0;\\
	\label{eqn:2.17}
	&(u(0),\theta(0))=(v_0,T_0).
\end{align}
Then the global well-posedness of 3D SPEs driven by fractional Brownian motion follows from the result in \cite{Zgl}.
\begin{thm}
	Let $Q_2\in L^2(\mho)$, $v_0\in\V_1,T_0\in\V_2, \T>0$. Assume the condition \eqref{eqn:2.10}  holds, then there exists a unique strong solution $(v,T)$ of the system \eqref{eqn:1.1}-\eqref{eqn:1.5}, or equivalently, $(u,\theta)$ of the system \eqref{eqn:2.12}-\eqref{eqn:2.17} on the interval $[0,\T]$ which is Lipschitz continuous with respect to the initial data and the noises in $\V$ and $C([0,\T];\V)$, respectively.
\end{thm}

%%%%%%%%%%%%%%%%%%%%%%%%
\section{Exponential stability of solutions}
\label{sec:3}
%%%%%%%%%%%%%%%%%%%%%%%

%%%%%%%%
\subsection{Steady state solution}
%%%%%%%
A stationary solution to \eqref{eqn:2.3} is $U^*=(v^*,T^*)$ satisfying
\begin{align*}
	&L_1v^*+(v^*\cdot \nabla)v^*-\left(\int_{-1}^z\nabla\cdot v^*(x,y,\xi,t)d\xi\right)\partial_zv^*+f(v^*)^{\perp}+\nabla p_s(x,y,z,t)-\int_{-1}^z\nabla T^*(x,y,\xi,t)d\xi=0,\notag\\
	&L_2T^*+(v^*\cdot \nabla )T^*-\left(\int_{-1}^z\nabla\cdot v^*(x,y,\xi,t)d\xi\right)\partial_z T^*=Q_2,\notag\\
\end{align*}
with boundary conditions \eqref{eqn:1.3}-\eqref{eqn:1.4} held. In abstract setting, the above equation can be rewritten as
\begin{equation}
\label{eqn:3.1}
	\nu AU^*+B(U^*)+R(U^*)=Q.
\end{equation}
The existence of steady state solutions follows from the result in \cite{M}.
\begin{thm}
	\label{thm:stablesolution}
	Suppose $Q_2\in L^2(\mho)$ and $\bar\nu>0$ is large enough, then \eqref{eqn:3.1}, together with boundary conditions \eqref{eqn:1.3}-\eqref{eqn:1.4}, has a unique solution $U^*=(v^*,T^*)$. Moreover, we have the following estimate
	\begin{equation}
		|AU^*|_2^2=|A_1v^*|_2^2+|A_2T^*|_2^2\leq K,
	\end{equation}
	where $K$ is a constant that depends on $Q,\nu$.
\end{thm}

%%%%%%%
\subsection{Exponential stability of steady state solutions}
%%%%%%%
In this section, we discuss the stability of steady state solutions. In the following, we first give the definition.
\begin{defn}
	We say the solution $U(t)$ to \eqref{eqn:2.3} converges to $U^*\in\H$ almost surely exponentially if there exists $\gamma>0$ such that
	\begin{equation}
		\lim_{t\to\infty}\frac{1}{t}\log|U(t)-U^*|_2\leq-\gamma.
	\end{equation}
	We say that $U^*$ is \textbf{almost surely exponentially stable} if any solution to \eqref{eqn:2.3} converges to $U^*$ almost surely exponentially with the same $\gamma>0.$
\end{defn}
The following Lemma is called uniform Gronwall lemma which will be used repeatedly in the proof of \autoref{prop:UH1}. One can refer to
Foias and Prodi \cite{FP} and Temam \cite{T} for a proof of \autoref{lemma:Gronwall}.
\begin{lemma}
	\label{lemma:Gronwall}
 Let $f,g$ and $h$ be three non-negative locally integrable functions on $(t_{0}, \infty)$ such that
\begin{eqnarray*}
 \frac{df}{dt}\leq gf +h,\ \ \ \forall\ t\geq t_{0},
\end{eqnarray*}
and
\begin{eqnarray*}
\int_{t}^{t+r}f(s)ds\leq a_{1},\ \ \int_{t}^{t+r}g(s)ds\leq a_{2},\ \ \int_{t}^{t+r}h(s)ds\leq a_{3},\ \ \forall\ t\geq t_{0},
\end{eqnarray*}
where $r, a_{1}, a_{2}, a_{3}$ are positive constants. Then
\begin{eqnarray*}
 f(t+r)\leq (\frac{a_{1}}{r}+a_{3} )e^{a_{2}},\ \ \ \forall\ t\geq t_{0}.
\end{eqnarray*}
\end{lemma}

\begin{prop}
	\label{prop:UH1}
	For any fixed $\omega$, there exists a random variable $C(\omega)$ taking values in $\mathbb{R}^{+}:=(0, \infty)$ such that
	\begin{align}
		\sup\limits_{t\in [0,\infty)}\|U(t)\|_1\leq C(\omega)<\infty,a.s..
	\end{align}
\end{prop}
\proof According to (5.106) in \cite{Zgl}, there exists positive constants $C$ and $\gamma_{1}$ such that
	\begin{align}
		\label{eqn:3.5}
&|u(t)|_2^2+|\theta(t)|_2^2\notag\\
\leq&(|u(0)|_2^2+|\theta(0)|_2^2)\exp\left\{\int_0^t-\gamma_1+C(\norm{Z_1}_2^2+\norm{Z_1}_2^4+\norm{Z_2}_3^2)ds\right\}\notag\\
&+\int_0^t\exp\left\{\int_s^t-\gamma_1+C(\norm{Z_1}_2^2+\norm{Z_1}_2^4+\norm{Z_2}_3^2)dx\right\}(|Q|_2^2+\norm{Z_1}_1^2)ds.
\end{align}
By \cite{MS}, we know that $(C_0(\R,\mathcal V),\mathcal B(C_0(\R,\mathcal V)), \P,\vartheta)$ is an ergodic metric dynamical system, it is shown in \cite{Zgl} that $Z_j(t)$ is adapted with respect to $\F_t:=\sigma(W_j^H(s),j=1,2, s\leq t)$. Therefore, applying the properties of ergodic metric dynamical system, there exists a $\beta$ big enough such that
$$
	\lim_{t\to\infty}\frac{C\int_0^t(\norm{Z_1}_2^2+\norm{Z_1}_2^4+\norm{Z_2}_3^2)ds}{t}=C\Ex[\norm{Z_1(0)}_1^2+\norm{Z_1(0)}_2^4+\norm{Z_2(0)}_3^2]<\frac{\gamma_1}{2}.
$$
Hence, for any fixed path $\omega$, there exists $T(\omega)$ big enough such that when $t>T(\omega)$,
$$\int_0^t(\norm{Z_1}_2^2+\norm{Z_1}_2^4+\norm{Z_2}_3^2)ds\leq \frac{\gamma_1}{2}t.$$
Now with \eqref{eqn:3.5}, for $t>T(\omega)$, we have
\begin{align}
\label{eqn:3.6}
	|U(t)|_2^2\notag=&|u(t)|_2^2+|\theta(t)|_2^2+|Z_1|_2^2+|Z_2|_2^2\notag\\
	\leq&(|u(0)|_2^2+|\theta(0)|_2^2)e^{-\gamma_1t/2}+\frac{1}{\gamma_1}e^{-\gamma_1t/2}\int_0^t(|Q|_2^2+\norm{Z_1}_1^2)ds\notag\\
	&+|Z_1|_2^2+|Z_2|_2^2.
\end{align}
By \autoref{prop:ZH3}, there exist $C_1(\omega), C_2(\omega)$ such that
$$\sup\limits_{t\in [0,\infty)}\norm{Z_1(t)}_1\leq C_1(\omega)<\infty \text{ and }\sup\limits_{t\in [0,\infty)}\norm{Z_2(t)}_1\leq C_2(\omega)<\infty.$$
Thus, back to \eqref{eqn:3.6},  by the continuity of $|U(t)|_{2}$ with respect to time $t,$ there exists $C(\omega),$
\begin{align}
\label{eqn:3.7}\sup\limits_{t\geq 0}|U(t)|_2^2 \leq C(\omega).
\end{align}
In view of (5.105) in \cite{Zgl}, (\ref{eqn:3.7}), \autoref{prop:ZH3} and \autoref{lemma:Gronwall}, we have
\begin{align}
\label{eqn:3.8} \int_{t}^{t+1}\|U(s)\|_{1}^{2}ds <C(\omega).
\end{align}
By the formula above (5.113) in \cite{Zgl}, (\ref{eqn:3.7}), (\ref{eqn:3.8}), \autoref{prop:ZH3} and \autoref{lemma:Gronwall}, we obtain
 \begin{align}
\label{eqn:3.9} \sup\limits_{t\in[0,\infty)}|\theta(t)|^2_{4} <C(\omega).
\end{align}
 Recall that $\bar u(x,y)=\int_{-1}^0u(x,y,z)dz$ and $\tilde u=u-\bar u$. By (5.122) in \cite{Zgl}, (\ref{eqn:3.7})-(\ref{eqn:3.9}), \autoref{prop:ZH3} and \autoref{lemma:Gronwall}, we obtain
  \begin{align}
\label{eqn:3.10} \sup\limits_{t\in[0,\infty)}|\tilde{u}(t)|^2_{4} <C(\omega).
\end{align}
 By virtue of (5.121) in \cite{Zgl}, (\ref{eqn:3.7})-(\ref{eqn:3.10}), \autoref{prop:ZH3} and \autoref{lemma:Gronwall}, we obtain
   \begin{align}
\label{eqn:3.11} \int_{t}^{t+1}\int_{\mho} \Big{(} |\nabla ( |\tilde{u}(s)|^{2})|^{2}+|\partial_{z}( |\tilde{u}(s)|^{2}) |^{2}ds
+\int_{t}^{t+1}\int_{\mho} |\tilde{u}|^{2}(| \nabla \tilde{u}(s)|^{2}+ | \partial_{z} \tilde{u}(s)|^{2} )ds <C(\omega).
\end{align}
 By (5.127) in \cite{Zgl}, (\ref{eqn:3.7})-(\ref{eqn:3.11}), \autoref{prop:ZH3} and \autoref{lemma:Gronwall}, we obtain
 \begin{align}
\label{eqn:3.12} \int_{t}^{t+1}|\Delta \bar{u} |_{2}^{2}ds<C(\omega)
 \end{align}
 and
  \begin{align}
\label{eqn:3.13}  \sup\limits_{t\in[0,\infty)}|\nabla \bar{u}(t) |_{2}^{2}ds<C(\omega).
 \end{align}
 By (5.137) in \cite{Zgl}, (\ref{eqn:3.7})-(\ref{eqn:3.13}), \autoref{prop:ZH3} and \autoref{lemma:Gronwall}, we obtain
  \begin{align}
\label{eqn:3.14} \int_{t}^{t+1}(|\nabla u_{z} |_{2}^{2}+ | u_{zz} |_{2}^{2}) ds<C(\omega)
 \end{align}
 and
  \begin{align}
\label{eqn:3.15}  \sup\limits_{t\in[0,\infty)}|u_{z}(t) |_{2}^{2}ds<C(\omega).
 \end{align}
 By (5.141) in \cite{Zgl}, (\ref{eqn:3.7})-(\ref{eqn:3.15}), \autoref{prop:ZH3} and \autoref{lemma:Gronwall}, we obtain
   \begin{align}
\label{eqn:3.16} \int_{t}^{t+1}|\Delta u (s) |_{2}^{2} ds<C(\omega)
 \end{align}
 and
  \begin{align}
\label{eqn:3.17}  \sup\limits_{t\in[0,\infty)}|\nabla u(t) |_{2}^2<C(\omega).
 \end{align}
 By (5.143) in \cite{Zgl}, (\ref{eqn:3.7})-(\ref{eqn:3.17}), \autoref{prop:ZH3} and \autoref{lemma:Gronwall}, we obtain
   \begin{align}
\label{eqn:3.18}  \sup\limits_{t\in[0,\infty)}(|\nabla \theta(t) |_{2}^2+|\theta_{z}|_{2}^2+\alpha|\theta(z=0)|_{2}^2)<C(\omega).
 \end{align}
 Finally, the result follows from (\ref{eqn:3.15}), (\ref{eqn:3.17}), (\ref{eqn:3.18}).
\qed

In the following, we give the main result of this paper.
\begin{thm}
	\label{thm:expstability}
	Suppose $Q_2\in L^2(\mho)$. We assume the conditions for $G_1,G_2$ stronger than \eqref{eqn:Gcon}, that is, there exist $M_1,M_2,\rho_1,\rho_2>0$ such that
\begin{equation}
\label{eqn:Gcon2}
|G_1(t)|+|G_1'(t)|\leq M_1e^{-\rho_1t},\  |G_2(t)|+|G_2'(t)|\leq M_2e^{-\rho_2t}.
\end{equation}
	 Let $U^*\in D(A)$ be the unique solution to \eqref{eqn:3.1} with boundary and initial conditions \eqref{eqn:1.3}-\eqref{eqn:1.5}, and let $\nu$ be large enough such that
	$$\min\{\nu_1,\mu_1,\nu_2,\mu_2\}>c_{0}\lambda_{1}|A U^{*}|_{2}^{2}+ \alpha_0\lambda_1,$$ the solution $U$ to \eqref{eqn:2.3} converges to the stationary solution $U^*$ almost surely exponentially. More precisely, there exists $0<\lambda<\rho_1\wedge \rho_2\wedge2\lambda_1^{-2}(\bar\nu-c_{0}\lambda_{1}|A U^{*}|_{2}^{2}- \alpha_0\lambda_1)$ such that
	\begin{equation}
		\lim_{t\to\infty}\frac{1}{t}\log|U(t)-U^*|_2\leq-\lambda/2.
	\end{equation}
\end{thm}
\proof Firstly, we have
\begin{align*}
	\frac{d}{dt}[e^{\lambda t}|U(t)-U^*|_2^2]=&\lambda e^{\lambda t}|U(t)-U^*|_2^2+2e^{\lambda t}\qv{\frac{d}{dt}U(t),U(t)-U^*}\\
	=&\lambda e^{\lambda t}|U(t)-U^*|_2^2-2e^{\lambda t}\nu\qv{AU(t),U(t)-U^*}\\
	&-2e^{\lambda t}\qv{B(U(t))+R(U(t)),U(t)-U^*}+2e^{\lambda t}\qv{Q(t),U(t)-U^*}\\
	&+2e^{\lambda t}\qv{G(t)d\dot{W}^H(t),U(t)-U^*}.
\end{align*}
By \eqref{eqn:3.1}, $U^*$ satisfies
$$e^{\lambda t}\qv{\nu AU^*+B(U^*)+R(U^*),U(t)-U^*}=e^{\lambda t}\qv{Q,U(t)-U^*},$$
we obtain that
\begin{align}
\label{eqn:3.21}
e^{\lambda t}|U(t)-U^*|_2^2
=&|U(0)-U^*|_2^2+\int_0^t\lambda e^{\lambda s}|U(s)-U^*|_2^2ds-2\int_0^te^{\lambda s}\nu\qv{A(U(s)-U^*),U(s)-U^*}ds\notag\\
&-2\int_0^te^{\lambda s}\qv{B(U(s))-B(U^*),U(s)-U^*}ds\notag\\
&-2\int_0^te^{\lambda s}\qv{R(U(s))-R(U^*),U(s)-U^*}ds\notag\\
&+2\int_0^te^{\lambda s}\qv{G(s)d{W}^H(s),U(s)-U^*}\notag\\
=&I_1+\cdots+I_6.
\end{align}
First by \eqref{eqn:L2H1}, we have
$$I_2\leq \int_0^t\lambda\lambda_1^{2} e^{\lambda s}\norm{U(s)-U^*}_1^2ds,\ I_3\leq-2\bar\nu\int_0^te^{\lambda s}\norm{U(s)-U^*}_1^2ds.$$
For the bilinear operator $B$, by \eqref{eqn:L2H1}, \autoref{lemma:2.1} and \autoref{thm:stablesolution}, we have
\begin{align*}
	I_4=&-2\int_0^te^{\lambda s}\qv{B(U(s),U(s)-U^*),U(s)-U^*}ds-2\int_0^te^{\lambda s}\qv{B(U(s)-U^*, U^*), U(s)-U^*}ds\\
	=&-2\int_0^te^{\lambda s}\qv{B(U(s)-U^*, U^*), U(s)-U^*}ds\\
	\leq&2c_0\lambda_1|AU^*|_2^2\int_0^te^{\lambda s}\norm{U(s)-U^*}_1^2ds
\end{align*}
By \eqref{eqn:L2H1} and \eqref{eqn:R},
$$I_5\leq 2\alpha_0\int_0^te^{\lambda s}\norm{U(s)-U^*}_1|U(s)-U^*|_2ds\leq 2\alpha_0\lambda_1\int_0^te^{\lambda s}\norm{U(s)-U^*}_1^2ds.$$
By the assumption, we see that for $\bar\nu>0$ large enough such that $\bar\nu>c_{0}\lambda_{1}|A U^{*}|_{2}^{2}+ \alpha_0\lambda_1$, then for small $\lambda$ so that
\begin{equation}
	\label{eqn:3.22}
	\lambda<2\lambda_1^{-2}(\bar\nu-c_{0}\lambda_{1}|A U^{*}|_{2}^{2}- \alpha_0\lambda_1),
\end{equation}
  we have
\begin{equation*}
	I_2+I_3+I_4+I_5<0.
\end{equation*}
By the definition of fractional Brownian motions in \eqref{eqn:2.6}, we get that
\begin{align*}
I_6=
	&2\int_0^te^{\lambda s}\qv{G_1(s)dW_1^H(s), v(s)-v^*}+2\int_0^te^{\lambda s}\qv{G_2(s)dW_1^H(s), T(s)-T^*}\\
	=&2\int_0^te^{\lambda s}\qv{\sum_{k=1}^{\infty}\sqrt{\lambda_{k,1}}e_{k,1}G_1(s)dB_k^H(s),v(s)-v^*}+2\int_0^te^{\lambda s}\qv{\sum_{k=1}^{\infty}\sqrt{\lambda_{k,2}}e_{k,2}G_2(s)dB_k^H(s),T(s)-T^*}\\
	=&2\sum_{k=1}^{\infty}\sqrt{\lambda_{k,1}}\int_0^t\qv{v(s)-v^*,e_{k,1}}e^{\lambda s}G_1(s)dB_k^H(s)+2\sum_{k=1}^{\infty}\sqrt{\lambda_{k,2}}\int_0^t\qv{T(s)-T^*,e_{k,2}}e^{\lambda s}G_2(s)dB_k^H(s).
\end{align*}
Now applying the inequality \eqref{eqn:2.5}, and by condition \eqref{eqn:Gcon},
\begin{align}
	&\bigg|\int_0^te^{\lambda s}\qv{G(s)dW^H(s), U(s)-U^*}\bigg|\notag\\
	\leq&\sum_{k=1}^\infty\sqrt{\lambda_{k,1}}\bigg|\sum_{j=1}^{[t]+1}\int_{j-1}^j\qv{v(s)-v^*,e_{k,1}}e^{\lambda s}G_1(s)dB_k^H(s)\bigg|\notag\\
	&+\sum_{k=1}^\infty\sqrt{\lambda_{k,2}}\bigg|\sum_{j=1}^{[t]+1}\int_{j-1}^j\qv{T(s)-T^*,e_{k,2}}e^{\lambda s}G_2(s)dB_k^H(s)\bigg|\notag\\
	\leq&\sum_{k=1}^\infty\sqrt{\lambda_{k,1}}\sum_{j=1}^{[t]+1}C_\alpha(B_k^H)\Big|_{j-1}^j \int_{j-1}^j\frac{|\qv{v(s)-v^*, e_{k,1}}|e^{-(\rho_1-\lambda )s}}{[s-(j-1)]^\alpha}ds\notag\\
	&+\sum_{k=1}^\infty\sqrt{\lambda_{k,2}}\sum_{j=1}^{[t]+1}C_\alpha(B_k^H)\Big|_{j-1}^j \int_{j-1}^j\frac{|\qv{T(s)-T^*, e_{k,2}}|e^{-(\rho_2-\lambda )s}}{[s-(j-1)]^\alpha}ds\notag\\
	&+\sum_{k=1}^\infty\sqrt{\lambda_{k,1}}\sum_{j=1}^{[t]+1}C_\alpha(B_k^H)\Big|_{j-1}^j\int_{j-1}^j\int_{j-1}^s\frac{|\qv{v(s)-v^*,e_{k,1}}e^{\lambda s}G_1(s)-\qv{v(x)-v^*,e_{k,1}}e^{\lambda x}G_1(x)|}{(s-x)^{1+\alpha}}dxds\notag\\
	&+\sum_{k=1}^\infty\sqrt{\lambda_{k,2}}\sum_{j=1}^{[t]+1}C_\alpha(B_k^H)\Big|_{j-1}^j\int_{j-1}^j\int_{j-1}^s\frac{|\qv{T(s)-T^*,e_{k,2}}e^{\lambda s}G_2(s)-\qv{T(x)-T^*,e_{k,2}}e^{\lambda x}G_2(x)|}{(s-x)^{1+\alpha}}dxds\notag\\
	=:&J_1+J_2+J_3+J_4.
\end{align}
Firstly, we have
\begin{align*}
	J_1\leq&\frac{1}{2}\sum_{k=1}^\infty\sqrt{\lambda_{k,1}}\sum_{j=1}^{[t]+1}C_\alpha(B_k^H)\Big|_{j-1}^j \bigg[\int_{j-1}^j|v(s)-v^*|_2^2e^{-(\rho_1-\lambda) s}ds+\int_{j-1}^j\frac{e^{-(\rho_1-\lambda) s}}{[s-(j-1)]^{2\alpha}}ds\bigg]\\
	\leq& \frac{1}{2}\sum_{k=1}^\infty\sqrt{\lambda_{k,1}}\sum_{j=1}^{[t]+1}C_\alpha(B_k^H)\Big|_{j-1}^j \bigg[\int_{j-1}^j|v(s)-v^*|_2^2e^{-(\rho_1-\lambda) s}ds+e^{-(\rho_1-\lambda) (j-1)}\frac{1}{1-2\alpha}\bigg]\\
	\leq&\frac{1}{2}\sum_{k=1}^{\infty}\sqrt{\lambda_{k,1}}\sum_{j=1}^{[t]+1}C_\alpha(B_k^H)\Big|_{j-1}^j \bigg[\int_{j-1}^j|v(s)-v^*|_2^2e^{-(\rho_1-\lambda) s}ds+\frac{e^{-(\rho_1-\lambda) (j-1)}}{1-2\alpha}\bigg].
\end{align*}
With similar discussion, we get
$$
	J_2\leq\frac{1}{2}\sum_{k=1}^{\infty}\sqrt{\lambda_{k,2}}\sum_{j=1}^{[t]+1}C_\alpha(B_k^H)\Big|_{j-1}^j \bigg[\int_{j-1}^j|T(s)-T^*|_2^2e^{-(\rho_2-\lambda) s}ds
	+\frac{e^{-(\rho_2-\lambda) (j-1)}}{1-2\alpha}\bigg].$$
Thus, by \autoref{prop:UH1}, for any $\omega$,
\begin{align}
\label{eqn:3.24}
	J_1+J_2\leq& \frac{1}{2}\sum_{k=1}^{\infty}(\sqrt{\lambda_{k,1}}+\sqrt{\lambda_{k,2}})\sum_{j=1}^{[t]+1}C_\alpha(B_k^H)\Big|_{j-1}^j \bigg[\int_{j-1}^j|U(s)-U^*|_2^2e^{-(\rho_1\wedge\rho_2-\lambda) s}ds
	+\frac{e^{-(\rho_1\wedge\rho_2-\lambda) (j-1)}}{1-2\alpha}\bigg]\notag\\
		\leq&\frac{1}{2}\sum_{k=1}^{\infty}(\sqrt{\lambda_{k,1}}+\sqrt{\lambda_{k,2}})\sum_{j=1}^{[t]+1}C_\alpha(B_k^H)\Big|_{j-1}^j[|U^*|_2^2+C(\omega)]e^{-(\rho_1\wedge\rho_2-\lambda)(j-1)}ds\notag\\
	&+\frac{1}{2}\sum_{k=1}^{\infty}(\sqrt{\lambda_{k,1}}+\sqrt{\lambda_{k,2}})\sum_{j=1}^{[t]+1}C_\alpha(B_k^H)\Big|_{j-1}^j\frac{e^{-(\rho_1\wedge\rho_2-\lambda) (j-1)}}{1-2\alpha}\notag\\
	=:&K_1+K_2.
\end{align}
By Lemma 7.4 and Lemma 7.5 in \cite{NR}, for any $0<\eps<H$,
\begin{align*}
	K_{1}+K_2\leq&\frac{1}{2}\sum_{k=1}^{\infty}(\sqrt{\lambda_{k,1}}+\sqrt{\lambda_{k,2}})\sum_{j=1}^{[t]+1}C_\alpha\eta_{\eps,j,j-1}\left(1+\frac{1}{H-\eps-1+\alpha}\right)[|U^{*}|_{2}^{2}+C(\omega)]e^{-(\rho_1\wedge\rho_2-\lambda)(j-1)}\\
	&+\frac{1}{2}\sum_{k=1}^{\infty}(\sqrt{\lambda_{k,1}}+\sqrt{\lambda_{k,2}})\sum_{j=1}^{[t]+1}C_\alpha\eta_{\eps,j,j-1}\left(1+\frac{1}{H-\eps-1+\alpha}\right)\frac{e^{-(\rho_1\wedge\rho_2-\lambda) (j-1)}}{1-2\alpha},
\end{align*}
where $\eta_{\eps,j,j-1}$ is a positive random variable such that $\Ex|\eta_{\eps,j,j-1}|^p<C_{\eps,p}$ that does not depend on $j$, for any $p\geq1$. Hence, for $\lambda<\rho_1\wedge\rho_2$,
\begin{align*}
	&\Ex\sum_{j=1}^{[t]+1}C_\alpha\eta_{\eps,j,j-1}\left(1+\frac{1}{H-\eps-1+\alpha}\right)e^{-(\rho_1\wedge\rho_2-\lambda)(j-1)}\\
	\leq&\sum_{j=1}^{\infty}C_\alpha C_{\eps,1}e^{-(\rho_1\wedge\rho_2-\lambda) (j-1)}<\infty,
\end{align*}
so we get $\sup\limits_{t\in [0,\infty)}[K_{1}+K_2]<\infty$, and back to \eqref{eqn:3.24}, we get $\sup\limits_{t\in [0,\infty)}[J_1+J_2]<\infty$.

Now for $J_3$, first by triangle inequality, and the assumption for $G_1$, we get
\begin{align*}
	J_3\leq&\sum_{k=1}^\infty\sqrt{\lambda_{k,1}}\sum_{j=1}^{[t]+1}C_\alpha(B_k^H)\Big|_{j-1}^j\int_{j-1}^j\int_{j-1}^s\frac{|\qv{v(s)-v(x),e_{k,1}}e^{\lambda x}G_1(x)|}{(s-x)^{1+\alpha}}dxds\\
	&+\sum_{k=1}^\infty\sqrt{\lambda_{k,1}}\sum_{j=1}^{[t]+1}C_\alpha(B_k^H)\Big|_{j-1}^j\int_{j-1}^j\int_{j-1}^s\frac{|\qv{v(s)-v^*,e_{k,1}}[e^{\lambda s}G_1(s)-e^{\lambda x}G_1(x)]|}{(s-x)^{1+\alpha}}dxds\\
	\leq&\sum_{k=1}^\infty\sqrt{\lambda_{k,1}}\sum_{j=1}^{[t]+1}C_\alpha(B_k^H)\Big|_{j-1}^j\int_{j-1}^j\int_{j-1}^s\frac{|\qv{v(s)-v(x),e_{k,1}}|e^{-(\rho_1-\lambda)x}}{(s-x)^{1+\alpha}}dxds\\
	&+\sum_{k=1}^\infty\sqrt{\lambda_{k,1}}\sum_{j=1}^{[t]+1}C_\alpha(B_k^H)\Big|_{j-1}^j\int_{j-1}^j\int_{j-1}^s\frac{|\qv{v(s)-v^*,e_{k,1}}|e^{-(\rho_1-\lambda)(j-1)}(s-x)}{(s-x)^{1+\alpha}}dxds\\
	\leq&\sum_{k=1}^\infty\sqrt{\lambda_{k,1}}\sum_{j=1}^{[t]+1}C_\alpha(B_k^H)\Big|_{j-1}^j\int_{j-1}^j\int_{j-1}^s\frac{|\qv{v(s)-v(x),e_{k,1}}|e^{-(\rho_1-\lambda)x}}{(s-x)^{1+\alpha}}dxds\\
	&+\sum_{k=1}^\infty\sqrt{\lambda_{k,1}}\sum_{j=1}^{[t]+1}e^{-(\rho_1-\lambda)(j-1)}C_\alpha(B_k^H)\Big|_{j-1}^j\int_{j-1}^j\int_{j-1}^s\frac{|\qv{v(s)-v^*,e_{k,1}}|}{(s-x)^{\alpha}}dxds,
\end{align*}
there is similar estimate for $J_4$, hence, altogether, by \autoref{prop:UH1},  one obtain that
\begin{align}
	J_3+J_4\leq&\sum_{k=1}^\infty(\sqrt{\lambda_{k,1}}+\sqrt{\lambda_{k,2}})\sum_{j=1}^{[t]+1}C_\alpha(B_k^H)\Big|_{j-1}^j\int_{j-1}^j\int_{j-1}^s\frac{|U(s)-U(x)|_2e^{-(\rho_1\wedge\rho_2-\lambda)x}}{(s-x)^{1+\alpha}}dxds\notag\\
	&+\sum_{k=1}^\infty(\sqrt{\lambda_{k,1}}+\sqrt{\lambda_{k,2}})\sum_{j=1}^{[t]+1}e^{-(\rho_1\wedge\rho_2-\lambda)(j-1)}C_\alpha(B_k^H)\Big|_{j-1}^j\int_{j-1}^j\int_{j-1}^s\frac{|U(s)-U^*|_2}{(s-x)^{\alpha}}dxds\notag\\
	=:&N_1+N_2.
\end{align}
Firstly, by \autoref{prop:UH1},
\begin{align*}
	N_2\leq
	&\sum_{k=1}^\infty(\sqrt{\lambda_{k,1}}+\sqrt{\lambda_{k,2}})\sum_{j=1}^{[t]+1}e^{-(\rho_1\wedge\rho_2-\lambda)(j-1)}C_\alpha(B_k^H)\Big|_{j-1}^j\int_{j-1}^j(C(\omega)+|U^*|_2)\frac{(s-j+1)^{1-\alpha}}{1-\alpha}ds\\
	\leq&\sum_{k=1}^\infty(\sqrt{\lambda_{k,1}}+\sqrt{\lambda_{k,2}})\sum_{j=1}^{\infty}\frac{(C(\omega)+|U^*|_2)e^{-(\rho_1\wedge\rho_2-\lambda)(j-1)}}{(2-\alpha)(1-\alpha)}C_\alpha(B_k^H)\Big|_{j-1}^j.
\end{align*}
With the same discussion as before, for $\lambda<\rho_1\wedge\rho_2$, we get $$\Ex \sum_{j=1}^{\infty}e^{-(\rho_1\wedge\rho_2-\lambda)(j-1)}C_\alpha(B_k^H)\Big|_{j-1}^j
<\infty,$$ so $\sup\limits_{t\in [0,\infty)}N_2<\infty,$ a.s..

Now we will estimate $N_{1}.$ Applying Theorem 4.2.1 in \cite{Z}, with $|B_k^H(t)|\leq t^2+c(\omega)$(see Lemma 2.6 of \cite{MS}), one can obtain
\begin{align*}
&\bigg|\qv{\int_{x}^{s}G_{1}\dot W_1^H,  e_{k,1}}\bigg|\\
&=\bigg|\langle \sum_{j=1}^\infty\sqrt{\lambda_{j,1}}e_{j,1}\int_x^sG_1(r)dB_j^H(r),    e_{k,1}    \rangle\bigg|=\sqrt{\lambda_{k,1}}\bigg|\int_x^sG_1(r)dB_k^H(r)\bigg|\\
&=\sqrt{\lambda_{k,1}}\bigg|G_1(s)B_k^H(s)-G_1(x)B_k^H(x)-\int_x^sG_1'(r)B_k^H(r)dr\bigg|\\
&\leq \sqrt{\lambda_{k,1}}[|G_{1}(s)B_k^{H}(s)- G_{1}(x)B_k^{H}(s) |+|G_{1}(x)B_k^{H}(s)- G_{1}(x)B_k^{H}(x) |+M_1(s^{2}+c(\omega))(s-x)]\\
&\leq \sqrt{\lambda_{k,1}}M_1[2(s-x)(s^{2}+c(\omega))+(s-x)^{H-\varepsilon}\eta_{\varepsilon, x,s }].
\end{align*}
Recall that $A_je_{j,k}=\gamma_{j,k}e_{j,k}$ for $j=1,2,$ then for $\eps<H-\alpha$ small enough,
\begin{align*}
	&|\qv{v(s)-v(x),e_{k,1}}|=\bigg|\int_x^s\qv{L_1v,e_{k,1}}+\int_x^s\qv{v\cdot\nabla v,e_{k,1}}+\int_x^s\qv{w\partial_zv,e_{k,1}}+\int_x^s\qv{G_1\dot W_1^H,e_{k,1}}\bigg|\\
	\leq&\gamma_{k,1}(|v|_2+|v|_2\norm{v}_1+{\norm{v}^{2}_1})(s-x)+\sqrt{\lambda_{k,1}}M_1[2(s-x)(s^{2}+c(\omega))+(s-x)^{H-\varepsilon}\eta_{\varepsilon, x,s }]\\
\leq&C\gamma_{k,1}(s-x)+C\sqrt{\lambda_{k,1}}[(s-x)(s^{2}+c(\omega))+(s-x)^{H-\varepsilon}\eta_{\varepsilon, x,s }].
\end{align*}
Similarly, we have
$$
	|\qv{T(s)-T(x),e_{k,2}}|\leq C\gamma_{k,2}(s-x)+C\sqrt{\lambda_{k,2}}[(s-x)(s^{2}+c(\omega))+(s-x)^{H-\varepsilon}\eta_{\varepsilon, x,s }].
$$
Hence, by the above estimates and \autoref{prop:UH1},
\begin{align*}
N_1\leq&C\sum_{k=1}^\infty(\sqrt{\lambda_{k,1}}\gamma_{k,1}+\sqrt{\lambda_{k,2}}\gamma_{k,2})\sum_{j=1}^\infty C_\alpha(B_k^H)\Big|_{j-1}^j\int_{j-1}^j\int_{j-1}^s\frac{e^{-(\rho_1\wedge\rho_2-\lambda)x}}{(s-x)^{\alpha}}dxds\\
&+C\sum_{k=1}^\infty(\lambda_{k,1}+\lambda_{k,2})\sum_{j=1}^{\infty}C_\alpha(B_k^H)\Big|_{j-1}^j\int_{j-1}^j\int_{j-1}^s\frac{e^{-(\rho_1\wedge\rho_2-\lambda)x}(s^2+c(\omega))}{(s-x)^{\alpha}}dxds\\
&+C\sum_{k=1}^\infty(\lambda_{k,1}+\lambda_{k,2})\sum_{j=1}^{\infty}C_\alpha(B_k^H)\Big|_{j-1}^j\eta_{\eps,j,j-1}\int_{j-1}^j\int_{j-1}^se^{-(\rho_1\wedge\rho_2-\lambda)x}(s-x)^{H-\eps-1-\alpha}dxds\\
\leq&C\sum_{k=1}^\infty(\sqrt{\lambda_{k,1}}\gamma_{k,1}+\sqrt{\lambda_{k,2}}\gamma_{k,2})\sum_{j=1}^\infty C_\alpha(B_k^H)\Big|_{j-1}^j\frac{e^{-(\rho_1\wedge \rho_2-\lambda)(j-1)}}{(1-\alpha)(2-\alpha)}\\
&+C\sum_{k=1}^\infty(\lambda_{k,1}+\lambda_{k,2})\sum_{j=1}^{\infty}C_\alpha(B_k^H)\Big|_{j-1}^je^{-\frac{1}{2}(\rho_1\wedge \rho_2-\lambda)(j-1)}\\
&+C\sum_{k=1}^\infty(\lambda_{k,1}+\lambda_{k,2})\frac{1}{(1+H-\eps-\alpha)(H-\eps-\alpha)}\sum_{j=1}^{\infty}C_\alpha(B_k^H)\Big|_{j-1}^je^{-(\rho_1\wedge\rho_2-\lambda)(j-1)}\eta_{\eps,j,j-1},
\end{align*}
where we used the boundedness $\sup\limits_{\{(x,s)\in [0,\infty)\times [0,\infty) \cap \{|x-s|\leq 1\}\}} e^{-\frac{1}{2}(\rho_1\wedge\rho_2-\lambda)x}(s^{2}+c(\omega)) <\infty$ in the last inequality.
Again by Lemma 7.4 and Lemma 7.5 in \cite{NR}, for $0<\eps<H$,
\begin{align*}
	N_1
	\leq&C\sum_{k=1}^\infty(\sqrt{\lambda_{k,1}}\gamma_{k,1}+\sqrt{\lambda_{k,2}}\gamma_{k,2})\sum_{j=1}^{\infty}\eta_{\eps,j,j-1}e^{-(\rho_1\wedge\rho_2-\lambda)(j-1)}\\
	&+C\sum_{k=1}^\infty(\lambda_{k,1}+\lambda_{k,2})\sum_{j=1}^{\infty}\eta_{\eps,j,j-1} e^{-\frac{1}{2}(\rho_1\wedge\rho_2-\lambda)(j-1)}  \\
	&+C\sum_{k=1}^\infty(\lambda_{k,1}+\lambda_{k,2})\sum_{j=1}^{\infty}e^{-(\rho_1\wedge\rho_2-\lambda)(j-1)}  \eta^2_{\eps,j,j-1},
\end{align*}
where $\eta_{\eps,j,j-1}$ is a positive random variable such that $\Ex|\eta_{\eps,j,j-1}|^p<C_{\eps,p}$ that does not depend on $j$, for any $p\geq 1$. Hence, there exists $C_2(\omega)$ such that
$$\Ex [\sup\limits_{t\in[0,\infty)}N_1]<C,$$
which implies
$$\sup\limits_{t\in[0,\infty)}N_1<C(\omega).$$
Back to \eqref{eqn:3.21}, we get that for any fixed $\omega$, there exists $C_1(\omega)$, and $\lambda<\rho_1\wedge\rho_2$ such that
\begin{equation*}
	e^{\lambda t}|U(t)-U^*|_2^2\leq C_1(\omega),
\end{equation*}
combining with \eqref{eqn:3.22}, this yields that for $\lambda<\rho_1\wedge \rho_2\wedge2\lambda_1^{-2}(\bar\nu-c_{0}\lambda_{1}|A U^{*}|_{2}^{2}- \alpha_0\lambda_1),$
\begin{equation}
	\lim_{t\to\infty}\frac{1}{t}\log|U(t)-U^*|_2^2\leq-\lambda.
\end{equation}
\qed


\begin{thebibliography}{99}


\bibitem {BCL}
Z. Brze$\acute{z}$niak, T. Caraballo b, J.A. Langa, Y. Li, G. {\L}ukaszewicz,
J. Real,   Random attractors for stochastic 2D-Navier Stokes
equations in some unbounded domains, J. Differential Equations, 255(2013), 3897--3919.


\bibitem {BLW}
Peter W. Bates, Kening Lu, Bixiang Wang, Random attractors for stochastic
reaction-Diffusion equations on unbounded domains, J. Differential Equations, 246(2009), 845--869.

\bibitem {CIN}
C. Cao, S. Ibrahim, K. Nakanishi, E.S. Titi, Finite-time blowup for the inviscid primitive
equations of oceanic and atmospheric dynamics, Comm. Math. Phys. 337(2015), 473--482.

\bibitem {CLT1}
C. Cao, J. Li, E.S. Titi, Global well-posedness of strong solutions to the 3D primitive equations
with horizontal eddy diffusivity, J. Differential Equations 257 (2014), 4108--4132.


\bibitem {CLT2}
C. Cao, J. Li, E.S. Titi, Local and global well-posedness of strong solutions to the 3D primitive
equations with vertical eddy diffusivity, Arch. Ration. Mech. Anal. 214 (2014), 35--76.

\bibitem {CLT3}
C. Cao, J. Li, E.S. Titi, Global well-posedness of the three-dimensional primitive equations
with only horizontal viscosity and diffusion, Communications on Pure and Applied Mathematics  Vol. LXIX(2016), 1492--1531.

\bibitem {CLTa1}
T. Caraballo, J. Langa, T. Taniguchi, The exponential behaviour and stabilizability of stochastic 2D--Navier--Stokes equations, J. Differential Equations, 179 (2002) 714--737.
\bibitem {CLTa2}
T. Caraballo, A. M. Marquez-Duran, J. Real, The asymptotic behaviour of a stochastic 3D LANS-$\alpha$ model, Appl. Math. Optim., 53(2006) 141--161.


\bibitem {CT1}
C. Cao, E.S. Titi, Global well-posedness of the three-dimensional viscous primitive equations of large scale ocean and atmosphere dynamics, Ann. of Math. 166(2007), 245--267.

\bibitem {CT2}
C. Cao, E.S. Titi, Global well-posedness of the 3D primitive equations with partial vertical
turbulence mixing heat diffusion, Comm. Math. Phys. 310 (2012), 537--568.

\bibitem {CDF}
H. Crauel, A. Debussche, F. Flandoli, Random attractors, J. Dynam. Differential Equations, 9(1997), 307--341.

\bibitem {CF}
H. Crauel, F. Flandoli, Attractors for random dynamical systems, Probab. Theory Related Fields, 100(1994), 365--393.

\bibitem {DGTZ}
A. Debussche, N. Glatt-Holtz, R. Temam, and M. Ziane, Global existence and regularity for the 3d stochastic primitive
equations of the ocean and atmosphere with multiplicative white noise, Nonlinearity, 25(2012), 2093.


\bibitem{DZr}
Z. Dong,  R. Zhang,
Markov selection and W-strong Feller for 3D stochastic primitive equations, Science China Mathematics, 60(2017), 1873--1900.


  \bibitem{DZZ1}  Z. Dong, J. Zhai, R. Zhang, Large deviation principles for 3D stochastic primitive equations,  J. Differential Equations 263(2017), 3110--3146.


\bibitem{DZZ2}
Z. Dong, J. Zhai, R. Zhang, Exponential mixing for 3D stochastic primitive equations of the large scale ocean. preprint. Available at arXiv: 1506.08514.




\bibitem {EPT}
B. Ewald , M. Petcu  and R. Temam , Stochastic solutions of the two-dimensional primitive equations
of the ocean and atmosphere with an additive noise,  Anal. Appl. (Singap.) 5(2007), 183--198.


\bibitem {FP}
C. Foias and G. Prodi, Sur le comportement global des solutions non-stationnaires des
\'equations de Navier-Stokes en dimension 2, Rend. Sem. Mat. Univ. Padova, 39(1967), 1--
34.


\bibitem{GGS}
Hongjun Gao, Maria J. Garrido-Atienza, and Bjorn Schmalfuss, Random attractors for stochastic evolution equations driven by fractional Brownian motion, SIAM J. MATH. ANAL., 46(2014), 2281--2309.


\bibitem {G}
A. E. Gill, Atmosphere-ocean dynamics, International Geophysics Series, Vol. 30, Academic Press,
San Diego, 1982.

\bibitem {GH}
B. Guo, D. Huang, 3d stochastic primitive equations of the large-scale ocean: global well-
posedness and attractors, Commun. Math. Phys. 286(2009), 697--723.

\bibitem {GHT1}
N. Glatt-Holtz, R. Temam , Cauchy convergence schemes for some nonlinear partial differential
equations, Appl. Anal. 90(2011), 85--102.

\bibitem {GHT2}
N. Glatt-Holtz,  R. Temam , Pathwise solutions of the 2-d stochastic primitive equations, Appl.
Math. Optim. 63(2011), 401--433.


\bibitem {GHZ}
N. Glatt-Holtz, M. Ziane,  The stochastic primitive equations in two space dimensions with
multiplicative noise, Discrete Contin. Dyn. Syst. Ser. B 10 (2008), 801--822.


\bibitem {GKVZ}
N. Glatt-Holtz, I. Kukavica, V. Vicol, and M. Ziane,
Existence and Regularity of Invariant Measures for the Three Dimensional Stochastic Primitive Equations, Journal of Mathematical Physics, 55(2014), 051504.


\bibitem{GS1}
Hongjun Gao, Chengfeng Sun, Well-posedness and large deviations for the stochastic primitive equations in two space dimensions, COMMUN. MATH. SCI.
 10(2012), 575--593.

\bibitem{GS2}
Hongjun Gao, Chengfeng Sun, Well-posedness of stochastic primitive equations with multiplicative noise in three dimensions, Disc. and Cont. Dyn. Sys. B, 21(2016), 3053--3073.

\bibitem{GS3}
Hongjun Gao, Chengfeng Sun, Hausdorff dimension of random attractor for stochastic Navier-Stokes-Voight equations and primitive equations. Dyn Partial Differ Equ, 7(2010), 307--326.



\bibitem {H1}
G. J. Haltiner, Numerical weather prediction, J.W. Wiley \& Sons,  New York, 1971.
\bibitem {H2}
G. J. Haltiner, R. T. Williams, Numerical prediction and dynamic meteorology, John Wiley \&
Sons, New York, 1980.

\bibitem{Ju} N. Ju, The global attractor for the solutions to the 3D viscous primitive equations, Discrete Contin. Dyn. Syst. 17 (1) (2007) 159--179.

\bibitem {KZ}
I. Kukavica  and M. Ziane, On the regularity of the primitive equations of the ocean, Nonlinearity
20(2007), 2739--2753.
\bibitem {LTW1}
J. L.  Lions,  R.  Temam and S.  Wang, New formulations of the primitive equations of atmosphere and applications, Nonlinearity 5(1992), 237--288.

\bibitem {LTW2}
J. L.  Lions,  R. Temam and S. Wang, On the equations of the large scale ocean, Nonlinearity  5(1992), 1007--1053.

\bibitem {LTW3}
J. L. Lions,  R. Temam and S. Wang, Models of the coupled atmosphere and ocean$(CAO I)$, Computational
Mechanics Advance 1 (1993), 1--54.

\bibitem {LTW4}
J. L. Lions,  R. Temam and S. Wang, Mathematical theory for the coupled atmosphere-ocean models
$(CAO III)$, J. Math. Pures Appl. 74 (1995), 105--163.


\bibitem{MN}
B. Maslowski and D. Nualart, Evolution equations driven by a fractional Brownian motion, Journal of Functional Analysis 202 (2003), 277--305.

\bibitem{MS} B. Maslowski, B. Schmalfuss, Random dynamical systems and stationary solutions of differential equations driven by the fractional Brownian motion, Stoch. Anal. Appl. 22 (6) (2004) 1577-1607.
	
\bibitem{M} T.T. Medjo, The primitive equations of the ocean with delays, Nonlinear Anal. Real World Appl. 10 (2) (2009) 779-797.
	
\bibitem{M1}
T.T. Medjo, The exponential behavior of the stochastic three-dimensional primitive equations with multiplicative noise, Nonlinear Anal. Real World Appl.,12(2011) 799--810.
	
\bibitem{MTV}
A.J. Majda, I. Timofeyev and E. Vanden-Eijinden,  A Mathematical Framework for Stochastic Climate Models, Commun.Pure Appl.Math., 54, 891--974, 2001.

\bibitem{NIVM} N. Glatt-Holtz, I. Kukavica, V. Vicol, and M. Ziane, Existence and regularity of invariant measures for the three dimensional stochastic primitive equations, Journal of Mathematical Physics, 55(2014), 051504.

	\bibitem{NR} D. Nualart, A. Rascanu, Differential equations driven by fractional Brownian motion, Collec. Math. 53 (2002) 55-81.



\bibitem{P}
T.N. Palmer, A nonlinear dynamical perspective on model error: a proposal for non-local stochastic-dynamic parametrizations in weather and climate prediction models,  Q.J.R.Meteorol.Soc., 127, 279--304, 2001.

	\bibitem{SKM} S. G. Samko, A. A. Kilbas, O. I. Marichev, Fractional Integrals and Derivatives: Theory and Applications, Gordon and Breach, London, 1993.

\bibitem{T}
R. Temam, "Infinite Dimensional Dynamical Systems in Mechanics and Physics," SpringVerlag, 1988, 2nd Edition, 1997.

	\bibitem{WZ} L. Wang, G. Zhou, The dynamical behavior of 3D stochastic primitive equations driven by fractional noise, \textit{preprint}.

	\bibitem{Z} M. Z\"ahle, Integration with respect to fractal functions and stochastic calculus, I, Probab. Theory and Related Fields 111 (1998) 333--374.
	


	\bibitem{Zgl} G. Zhou, Random attractor of the 3D viscous primitive equations driven by fractional noises, J. Differential Equations, 266(2019), 7569--7637.
\end{thebibliography}
\end{document}